%% file: pltinf.tex
\documentclass[a4paper,12pt]{article}

\usepackage{amsmath}
\usepackage{amsfonts}
\usepackage{amssymb}
\usepackage{amsthm}
\usepackage{graphics}
\usepackage{amscd}
\usepackage{graphicx,epsfig}

\DeclareMathOperator{\Div}{div}
\DeclareMathOperator{\im}{im}
\DeclareMathOperator{\Gra}{Graph}

\renewcommand{\epsilon}{\varepsilon}
\newcommand{\boF}{\mathcal{F}}

\newcommand{\boM}{\mathcal{M}}

\newcommand{\boG}{\mathcal{G}}
\newcommand{\boB}{\mathcal{B}}

\newcommand{\boP}{\mathcal{P}}

\newcommand{\R}{\mathbb{R}}

\newcommand{\C}{\mathbb{C}}
\renewcommand{\L}{\mathbb{L}}

\newcommand{\N}{\mathbb{N}}

\newcommand{\der}[2]{\dfrac{\partial #1}{\partial #2}}
\newcommand{\dd}{\mathrm{d}}

\newtheorem{defn}{Definition}
\newtheorem{thm}{Theorem}

\newtheorem{prop}{Proposition}
\newtheorem{lem}{Lemma}

\renewcommand{\phi}{\varphi}
\newcommand{\dis}{\displaystyle}

\newcounter{remark}

\newenvironment{rem}[1]{\refstepcounter{remark}\label{#1}
  \noindent\textbf{Remark \ref{#1}.}}{}

\usepackage{color}

\newcommand{\Ome}{\Omega}

\title{The Dirichlet problem for minimal surfaces equation and Plateau problem
  at infinity}
\author{Laurent MAZET
 \footnote{Laboratoire Emile Picard (UMR 5580), Universit\'e
    Paul Sabatier, 118, route de Narbonne, 31062 Toulouse, France; e-mail
    mazet@picard.ups-tlse.fr}} 
\date{}

\begin{document}
\maketitle

\begin{abstract}
In this paper, we shall study the Dirichlet problem for the minimal
surfaces equation. We prove some results about the boundary behaviour
of a solution of this problem. We describe the behaviour of a
non-converging sequence of solutions in term of lines of divergence in
the domain. Using this second result, we build some solutions of the
Dirichlet problem on unbounded domain. We then give a new proof of the
result of C.~Cos\'\i n and A.~Ros concerning the Plateau problem at
infinity for horizontal ends.
\end{abstract}



\setcounter{section}{-1}

\section{Introduction}

\input pltinf0.tex

\section{The Dirichlet problem on multi-domains} 

In this section, we shall give a generalization of the results of
H. Jenkins and J. Serrin \cite{JS} for the Dirichlet problem on bounded
domain. First we have to generalize the notion of domain of
$\R^2$. Let us consider a pair $(\Omega,\phi)$ where $\Omega$ is a 
simply-connected 2-dimensional complete flat manifold with piecewise
smooth boundary and $\phi:\Omega\longrightarrow\R^2$ is a local
isometry. The map $\phi$ is called the \emph{developing map} and the
points where the boundary $\partial\Omega$ are not smooth are called
\emph{vertices}.
\begin{defn}\label{multdom}
A   pair  $(\Omega,\phi)$,  where   $\Omega$  is   a  simply-connected
2-dimensional  complete flat manifold  with piecewise  smooth boundary
and  $\phi:\Omega\longrightarrow\R^2$  is   a  local  isometry,  is  a
\emph{multi-domain} if each connected  component of the smooth part of
$\partial \Omega$ is a convex arc.
\end{defn}

If $(\Ome,\phi)$ is as above and a part of $\partial\Omega$ is linear
then we add two vertices at 
the end points of this linear part and we call this new part an
\emph{edge}. 

Let $(\Omega,\phi)$ be a multi-domain, if $u$ is a smooth function on
$\Omega$ we shall call \emph{graph} of $u$ the surface in $\R^3$ given
by
$\{(\phi(x),u(x))\}_{x\in\Omega}$. If $u$ is a solution of the minimal
surfaces equation \eqref{MSE}, the graph of $u$ is a minimal surface
of $\R^3$. The Dirichlet problem on multi-domain consists in the
determination of a 
function $u$ satisfying the equation \eqref{MSE} on $\Omega$ and
taking on assignated values on the boundary of $\Omega$.

As in the case of a domain in $\R^2$, if $u$ is a solution of
\eqref{MSE} on $\Omega$, we can define a differential form $\dd \Psi_u$
on $\Omega$ which corresponds to the differential of the third
coordinate of the conjugate surface of the graph of $u$. In using the
charts given by 
the developing map $\phi$, we have $\dd\Psi_u= \dis\frac{p}{W}\dd y
-\frac{q}{W}\dd x$. $\dd \Psi_u$ is a closed form by \eqref{MSE} and, since 
$\Omega$ is simply connected, we can define a function $\Psi_u$ on
$\Omega$ which is $1$-Lipschitz continuous, we call this function the
\emph{conjugate function} to $u$. One important result concerning
$\dd \Psi_u$ is the following lemma.

\begin{lem}{\label{calc}}
Let $\Omega$ be a domain bounded in part by a straight segment $T$, oriented
such that the right hand normal to $T$  is the outer normal to $\Ome$. Let $u$
be a solution of \eqref{MSE} in $\Ome$ which assumes the boundary value $+\infty$
on $T$. Then
$$\int_T \dd\Psi_u=|T|$$
\end{lem}
This is Lemma 4 in \cite{JS}. For other properties of $\Psi_u$ and $\dd
\Psi_u$, we refer to this paper. 

When $\Omega$ is compact there is a finite number of connected
components of the smooth part of $\partial\Omega$; let us call them
$C_1,\dots,C_n$. When the data on the boundary is bounded, we have
this result:

\begin{thm}
Let $(\Omega,\phi)$ be a compact multi-domain with boundary arcs
$C_1$,\dots, $C_n$ and let $u_1,\dots,u_n$
be bounded continuous functions respectively on $C_1$,\dots, $C_n$. Then
there exists a unique solution $u$ of the minimal surfaces equation on
$\Omega$ such that $u_{|C_i}=u_i$.
\end{thm}   

\begin{proof}
The proof of the uniqueness is a particular case of the proof of Theorem
\ref{thmjs}, so we make it later.

The existence of the solution on multi-domain is due to a Perron
process, let us recall some elements of this method. If $v$ is a
continuous function on $\Omega$ and $D$ is a disk in $\Omega$, we note
by $u_{v,D}$ the solution of \eqref{MSE} in $D$ which takes the value
$v$ on $\partial D$. We also note $M_D[v]$ the continuous function
which coincides with $v$ on $\Omega\backslash D$ and $u_{v,D}$ on
$D$. Let $u_1,\dots,u_n$ be the data on the boundary of $\Omega$; we
say that $v$ is a sub-solution of the Dirichlet problem if $v\le u_i$
on $C_i$ and $v\le M_D[v]$ for all disks $D$ in $\Omega$. Since the
$u_i$ are bounded by a constant $M$, the class $\mathcal{F}$ of all
sub-solutions is non-empty: the constant function $-M$ is in; besides,
each sub-solution $v$ verifies $v\le M$. So we can define a function
$u$ by:

\begin{equation}
\forall P\in\Omega \quad\quad u(P)=\sup_{v\in\mathcal{F}}v(P)
\end{equation}

By standard argument, we can show that u is a solution of
\eqref{MSE}. Since in our definition of multi-domain we suppose that
the boundary is locally convex, there exist barrier functions on the
boundary (they are constructed in using the Scherk surface). So we can
insure that $u$ takes the value $u_i$ on $C_i$. For more details on Perron
process, we can refer to the book of D.~Gilbarg and N.S.~Trudinger \cite{GT}
or the one of R.~Courant and D.~Hilbert \cite{CH} which illustrate
this method for the classical Laplacian Dirichlet problem, there is
also the book of J. C. C. Nitsche \cite{Ni2} which studies the case of
the minimal surfaces equation. 
\end{proof}

The work of H. Jenkins and J. Serrin is to allow infinite data on the
boundary. By the Straight Line Lemma \cite{JS}, we know that infinite
data can only be allowed on linear parts of the boundary.
\begin{defn}
Let $(\Omega,\phi)$ be a multi-domain, a \emph{polygonal domain} $\mathcal{P}$
of $\Omega$ is a connected compact subset of $\Omega$ such that
$(\mathcal{P},\phi)$ is a multi-domain, the boundary of $\mathcal{P}$
is only composed of edges and the vertices of $\mathcal{P}$ are drawn
from the vertices of $\Omega$.
\end{defn}
We want to solve the Dirichlet problem with infinite data so let us
call $A_1,\dots,A_k$ and $B_1, \dots, B_l$ the edges of $\Omega$ such
that we assign the value $+\infty$ on $A_j$ and $-\infty$ on
$B_j$. We call $C_1,\dots,C_n$ the remaining arcs on which we
assign continuous data. 

Let $\mathcal{P}$ be a polygonal domain of $\Omega$. We note,
respectively,  $\alpha$ and $\beta$ the total length of the edges
$A_j$ and the one of the edges $B_j$ which belong to the
boundary of $\mathcal{P}$ and we note $\gamma$ the perimeter of
$\mathcal{P}$. We then have the following generalization of the result
of H. Jenkins and J. Serrin.

\begin{thm}{\label{thmjs}}
Let $(\Omega,\phi)$ be a compact multi-domain with the families
$\{A_j\}$, $\{B_j\}$ and $\{C_j\}$ as above.

If the familly $\{C_j\}$ is non-empty, then there exists a solution of
the minimal surface equation in $\Omega$ which assumes the value
$+\infty$ on each $A_j$, the value $-\infty$ on each $B_j$ and
arbitrarily assignated continuous data on each $C_j$, if and only if

\begin{equation}{\label{cond}}
2\alpha<\gamma\quad \textrm{and} \quad 2\beta<\gamma
\tag{$*$}
\end{equation}
for each polygonal domain $\mathcal{P}$ of $\Omega$. If a solution
exists, it is unique.

If the familly $\{C_j\}$ is empty, then a solution exists, if and only
if

$$\alpha=\beta$$
when $\mathcal{P}$ coincides with $\Omega$ and \eqref{cond} holds for
all other polygonal domains of $\Omega$. In this case, if a solution
exists, it is unique up to an additive constant. 
\end{thm}

\begin{proof}
To prove the existence of a solution, we can use the same arguments
than H. Jenkins and J. Serrin, so we refer to \cite{JS}.

The proof of the uniqueness in \cite{JS} works also but we give another
proof which we can apply in other situations. Let $u_1$ and $u_2$ be
different solutions of \eqref{MSE} with the same data on the
boundary. In the case where the familly $(C_j)$ is empty, we suppose that
$u_1-u_2$ is not constant; besides, in considering $u_i-u_i(P)$ (where
$P\in \Omega$), we can assume that $\{u_1<u_2\}$ and $\{u_1>u_2\}$ are
non-empty. In choosing sufficiently small $\epsilon>0$, we have
$\Omega_\epsilon= \{u_1-u_2>\epsilon\} \neq\emptyset$, besides the choise of
$\epsilon$ is such that $\partial\Omega_\epsilon$ is regular. We note
$\dd\tilde{\Psi}=\dd\Psi_{u_1}-\dd\Psi_{u_2}$, since $\dd\tilde{\Psi}$
is closed, we have $\dis \int_{\partial\Omega_\epsilon}
\dd\tilde{\Psi}= 0$. Because $u_1$ and $u_2$ have the same data on the
boundary, $\partial\Omega_{\epsilon}$ does not intersect $\dis\cup_j
C_j$ so $\partial\Omega_{\epsilon}$ is composed of three parts: one is
included in 
$\dis\cup_j A_j\bigcup \cup_j B_j$ on which $\dd\tilde{\Psi}=0$ (this is a
consequence of Lemma \ref{calc}), one is included in $\Omega$ and a
last part which is composed of 
some vertices of $\Omega$ but its contribution to the integral is
zero. On the second part, let us call it
$\widetilde{\partial\Omega_{\epsilon}}$, $\nabla u_1-\nabla u_2$ points in
$\Omega_{\epsilon}$, this part is then oriented by the non-direct
normal to  $\nabla u_1-\nabla u_2$ so, by Lemma $2$ of P. Collin
and R. Krust in \cite{CK}, $\dis \int_{\widetilde{\partial\Omega_{\epsilon}}}
\dd\tilde{\Psi}< 0$; this gives us a contradiction.
\end{proof}

\section{A result of regularity at the vertices}

\input{pltinf2.tex}

\section{Convergence and divergence of sequence of solutions of \eqref{MSE} 
{\label{seqmse}}} 

In this section we shall explain what we can say when we have a
sequence $(u_n)$ of solutions of \eqref{MSE} about its convergence: can
we make converge a subsequence by some compactness result? What are the
different ways of divergence? In 
\cite{JS}, it is shown that for a monotone sequence, it appears lines
which separate domains of convergence and domains of divergence (this works
only for subsequence). We
shall show that such lines always appear (Theorem \ref{linofdiv2}).

First, we have to determine
the domain on which we can make converge a sequence, since each
surface is a graph, if we want the limit to be a graph, the normal
to the surface needs to stay close by the vertical unit vector and then
$W_n$ have to be bounded. We have then the following lemma.

\begin{lem}{\label{bound}}
Let $\Ome$ be a domain and $(u_n)$ a sequence of solutions of \eqref{MSE}
on $\Ome$. Let $P\in \Ome$; we suppose that $W_n(P)$ is bounded by a
constant $M$; then there exists $R>0$ which depends only of $M$ and the
distance of $P$ to $\partial \Ome$ such that $W_n$ is bounded by $2M$ on
the disk of center $P$ and radius $R$.
\end{lem}

\begin{proof}
We fix an index $n$. We know (see \cite{Ni}) that there exists a
constant $c$ such that if $u$ is a solution of \eqref{MSE} on the disk
$\{(x,y)|\ x^2+y^2<R^2\}$ we have:
\begin{equation}
r^2(0)+2s^2(0)+t^2(0)\le\frac{c}{R^2}W^4(0)
\end{equation}
Let $R$ be such that $2R= d(P,\partial \Ome)$ then, for all $Q$ in
$D(P,R)$, the above equation gives $\dis
r^2(Q)+2s^2(Q)+t^2(Q)\le\frac{c}{R^2}W^4(Q)$. We 
have $\dis\nabla W=(\frac{rp+sq}{W},\frac{sp+tq}{W})$, so, in $D(P,R)$, we
have $||\nabla W||\le \widetilde{C}W^2$ with $\widetilde{C}$ which
depends only of $R$. Let $z$ be the function such that $z(0)=M$ and
$z'=\widetilde{C}z^2$, $z$ is defined on $\dis[0,
\frac{1}{M\widetilde{C}}[$ by: 
\begin{equation}
\frac{1}{M}-\frac{1}{z}=\widetilde{C}r
\end{equation}
Because of our estimate on $||\nabla W||$, we have, in polar
coordinates, $W(r,\theta)\le z(r)$. Then $W$ is bounded by $2M$ on
$\dis D(P,\min(R,\frac{1}{2M\widetilde{C}}))$. 
\end{proof}

Let $(u_n)$ be a sequence of solutions of \eqref{MSE} on a domain
$\Ome$. We then define $\boB(u_n)$ as the set of the point $Q\in \Ome$ such
that $(W_n(Q))$ is bounded. Lemma \ref{bound} says us that $\boB(u_n)$
is an open set and that $W_n$ is uniformly bounded on each compact
inclued in $\boB(u_n)$. Then if $D$ is a connected component of
$\boB(u_n)$ and $P\in D$ there exists an extraction $\theta$ such
that $u_{\theta(n)}-u_{\theta(n)}(P)$ converges uniformly on each compact
  of $D$ to a solution $u$ of \eqref{MSE}; here, we use some
  classical compactness results (see \cite{Ni}). This proves that the
  divergence of the sequence is due to the behaviour of the sequence
  over $\Ome\backslash \boB(u_n)$.

If $P\in \Ome\backslash \boB(u_n)$, there exists a subsequence $u_{n'}$
such that $W_{n'}(P)\longrightarrow +\infty$. As the normal $N_n$ to
the graph at $(P,u_n(P))$ is given by: 
\begin{equation}{\label{norm}}
N_n(P)=\left(\frac{p_n}{W_n}(P),\frac{q_n}{W_n}(P),
  -\frac{1}{W_n}(P)\right)
\end{equation}
we can suppose that $N_{n'}(P)$ converges to an horizontal unit
vector. The following proposition describes what locally happens.
\begin{prop}{\label{linofdiv1}}
Let $r$ be positive. Let $(u_n)$ be a sequence of solutions of
\eqref{MSE} on the disk $D(0,r)$. We suppose that $N_n(0)$ converges to
$(1,0,0)$. Let $\alpha\in]0,1[$, then there exists an extraction
$\theta$ such that $N_{\theta(n)}$ converges to $(1,0,0)$ at almost
every point of $\{0\}\times [-\alpha r,\alpha r]$.
\end{prop}

\begin{proof}
Let $n\in\N$, we know (see \cite{Os} and \cite{JS}) that there exists
$\Phi_n:(x,y) \mapsto (\xi, \eta)$ with $\Phi_n(0,0)=(0,0)$ and:
\begin{gather}
\dd \xi= \left(1+\frac{1+p_n^2}{W_n}\right)\dd x+
\frac{p_nq_n}{W_n}\dd y\\
\dd \eta= \frac{p_nq_n}{W_n}\dd x+
\left(1+\frac{1+q_n^2}{W_m}\right)\dd y
\end{gather}

We know that $\Phi_n$ increases distance so it is bijective on its
image. This image contains the disk of center $(0,0)$ and radius
$r$. Besides, we know that $(\xi,\eta)$ are conformal parameters for the
 graph of $u_n$. On the $\xi\eta$ disk $D(0,r)$ we then have the Gauss
map $g_n(\xi+i\eta)$ which corresponds to the stereographic projection
of $N_n$; $g_n$ is holomorphic. We note $z_n=g_n(0)$, by hypothesis we
have $z_n\longrightarrow 1$. We note $z=\xi+i\eta$; by our choice of normal
$g_n:D(0,r)\longrightarrow D(0,1)$, then there exists
$h_n:D(0,r)\longrightarrow D(0,1)$ holomorphic with $h_n(0)=0$ such
that:
\begin{equation}
g_n(z)=\frac{h_n(z)+z_n}{1+\overline{z_n}h_n(z)}
\end{equation}

Since $z_n\longrightarrow 1$, the sequence of holomorphic functions
$\dis z\longmapsto\frac{z+z_n}{1+\overline{z_n}z}$ converges simply to
$1$ on $D(0,1)$ and uniformly on the disk $D(0,\alpha)$ for all
$\alpha<1$. But by Schwarz Lemma, we have, for all $n\in\N$,
$h_n\big(D(0,\alpha r)\big)\subset D(0,\alpha)$, we then have uniform
convergence of $g_n$ to $1$ on $D(0,\alpha r)$. In using \eqref{norm}, this
proves that for every $\epsilon$, if $n$ is big enough, we can say that:
$\dis\frac{p_n}{W_n}\ge 1-\epsilon$ and $ \dis \frac{|q_n|}{W_n}\le \epsilon$
in $\Phi_n^{-1}\left(D(0,\alpha r)\right)$. So to conclude, we need to
understand the shape of $\Phi_n^{-1}\left(D(0,\alpha r)\right)$; we
shall see that these sets are concentrating along the segment
$\{0\}\times[-\alpha r,\alpha r]$.

For all $n\in \N$, we consider, in the $\xi\eta$ disk, the path
$\gamma_n:]-r,r[ 
\longrightarrow D(0,r)$ defined by $\gamma_n(0)=0$ and
$\gamma_n'=\frac{\nabla y_n}{||\nabla y_n||}$ where $y_n$ is the
second coordinate of $\Phi_n^{-1}$. We have (see \cite{Os}):
\begin{equation}
\nabla y_n=\left(-\frac{p_nq_n}{J_nW_n},\frac{W_n+1+p_n^2}{J_nW_n}\right)
\end{equation}
where $J_n=\det(\dd\Phi_n)=\dis W_n+2+\frac{1}{W_n}$.Because \eqref{norm},
$\nabla y_n$ converges 
uniformly to $(0,1)$ on the disk $D(0,\alpha r)$ for all $\alpha<1$.

Let $\alpha\in]0,1[$, we note $A_\alpha$ and $B_\alpha$ the points in
the $xy$-disk $D(0,r)$ of respective coordinates $(0,-\alpha r)$ and
$(0,\alpha r)$. In the following, we prove that $\dis \int_{[A_\alpha,
  B_\alpha]} \dd\Psi_{u_n}\longrightarrow 2\alpha r$.

Let $\widetilde{\alpha}>\alpha$, then for $n$ big enough, we have
$||\nabla y_n||>\dis\frac{\alpha}{\widetilde{\alpha}}$ in $D(0,
\widetilde{\alpha}r)$. Because, for $t\in[-\widetilde{\alpha}r,
\widetilde{\alpha}r]$, $\gamma_n(t)\in D(0,\widetilde{\alpha}r)$,
there exists $\widetilde{\alpha}r\le t_0^n< t_1^n\le \widetilde{
  \alpha}r$ such that $y_n( \gamma_n( t_0^n))= -\alpha r$ and $y_n(
\gamma_n( t_1^n))= \alpha r$. Along $[t_0^n,t_1^n]$, $y_n\circ
\gamma_n$ increases strictly from $-\alpha r$ to $\alpha r$, then the
path $\Gamma_n= \Phi_n^{-1}\circ\gamma_n$ on $[t_0^n,t_1^n]$ can be
parametrized by $y\in[-\alpha r,\alpha r]$: we have a function $f_n$
on $[-\alpha r,\alpha r]$ such that for $t\in[t_0^n,t_1^n]$ $x_n
(\gamma_n(t))= f_n\big(y_n(\gamma_n(t))\big)$. We have:
\begin{equation}
\begin{split}
|x_n\circ \gamma_n(t)|&\le
\left|\int_0^t|(x_n\circ\gamma_n)'(u)|\dd u \right|\\
&\le\left|\int_0^t||\nabla x_n||(\gamma_n(u))\dd u \right| =
\left|\int_0^t \left( \frac{1+q_n^2}{(1+W_n)^2} \right)^\frac{1}{2}
  (\Gamma_n(u))\dd u\right| 
\end{split}
\end{equation}
(for the last equality see \cite{Os}), then for $n$ big enough
$|x_n\circ\gamma_n|$ on $[t_0^n,t_1^n]$ can be bounded by a constant uniformly
small; this is due to the fact that $ \dis \frac{|q_n|}{W_n}\le \epsilon$ for
big $n$. We then have proved that the path $\Gamma_n$ is close by the
segment $[A_\alpha, B_\alpha]$ for big $n$.

Let us now calculate $\dis\int_{\Gamma_n} \dd\Psi_{u_n}$. We have: 
\begin{equation}
\int_{\Gamma_n}\frac{p_n}{W_n}\dd y-\frac{q_n}{W_n}\dd x=\int_{-\alpha
  r}^{\alpha r}\left(\frac{p_n}{W_n}(f_n(y),y)-
  \frac{q_n}{W_n}(f_n(y),y) f_n'(y)\right)\dd y 
\end{equation}
We have $|f_n'\big(y_n(\gamma_n(t))\big)|\le \dis\frac{||\nabla
  x_n||}{||\nabla y_n||}(\gamma_n(t))\longrightarrow 0$, the
  convergence is uniform so $f_n'$ tends uniformly to $0$ on $[-\alpha
  r,\alpha r]$. This proves that $\dis \int_{\Gamma_n}
  \dd\Psi_{u_n}\longrightarrow 2\alpha r$

Let us consider the path $\widetilde{\Gamma_n}$ which consists on the
segment $[A_\alpha,B_\alpha]$, then the segment
$[B_\alpha,\Gamma(t_1^n)]$, then $-\Gamma_n$, then, finally, the segment
$[\Gamma_n(t_0^n),A_\alpha]$. Let $\epsilon>0$, for $n$ big enough, we
can suppose that $\dis\int_{\Gamma_n} \dd\Psi_{u_n}\ge 2\alpha
r-\epsilon$ and $\dis \left( \frac{1+q_n^2}{(1+W_n)^2}
\right)^\frac{1}{2} \le \epsilon$ on
$\Phi_n^{-1}\big(D(0,\widetilde{\alpha} r)\big)$. As $\dd \Psi_{u_n}$
is closed $\dis \int_{\widetilde{\Gamma_n}} \dd\Psi_{u_n}=0$; we then
have:
\begin{equation}
\begin{split}
2\alpha r\ge \int_{[A_\alpha,B_\alpha]}\dd\Psi_{u_n}
&=-\int_{[B_\alpha,\Gamma(t_1^n)]}\dd \Psi_{u_n}
-\int_{-\Gamma_n}\dd \Psi_{u_n}
-\int_{[\Gamma_n(t_0^n),A_\alpha]}\dd \Psi_{u_n}\\
&\ge \int_{\Gamma_n} \dd\Psi_{u_n} -|x_n\big(\gamma(t_1^n)\big)|
-|x_n\big(\gamma(t_0^n)\big)|\\
&\ge 2\alpha r-\epsilon -\epsilon |t_1^n|-\epsilon |t_0^n|\\
&\ge 2\alpha r-\epsilon (1+2\widetilde{\alpha}r)
\end{split}
\end{equation}

This proves that $\dis\int_{[A_\alpha,B_\alpha]}\dd
\Psi_{u_n}\longrightarrow 2\alpha r$ as $n$ tends to $+\infty$. We have
$\dis\int_{[A_\alpha,B_\alpha]}\dd \Psi_{u_n}=\int_{[A_\alpha,
  B_\alpha]} \frac{p_n}{W_n}\dd y$. Because $\dis\frac{p_n}{W_n}\le
1$, the preceding equalities prove that $\dis\frac{p_n}{W_n}$
converges to $1$ in $\L^1([A_\alpha, 
B_\alpha])$. Then there exists an extraction $\theta$ such that $\dis
\frac{p_{\theta(n)}}{W_{\theta(n)}}$ converges simply to $1$ at almost
every point in $[A_\alpha, B_\alpha]$; thus the proposition is
proved. 
\end{proof}

This proposition gives us a local result and we have the following
global result

\begin{thm}\label{linofdiv2}
Let $(\Omega,\phi)$ be a multi-domain. Let $(u_n)$ be a sequence of
solutions of \eqref{MSE} on $\Omega$. Let $P\in\Omega$ and $N$ a unit
tangent vector at $P$, we call $D$ the geodesic of $\Omega$ passing at
$P$ and normal to $N$. If the sequence $(N_n(P))$ converges to $N$,
then $N_n(Q)$ converges to $N$ at every point of $D$.
\end{thm}

As $\Omega$ is locally isometric to $\R^2$, we have allowed us to call
$N$ the parallel transport of $N$ along $D$.

\begin{proof}
We first get a parametrization of $D$ by arc-length with $P$ as
origin-point; then $D$ is parametrized by $]a,b[$, $]-\infty,b[,
]a,+\infty[$ or $]-\infty,+\infty[$, we shall suppose that we are in
the case $]a,+\infty[$ (the other cases are similar). We then consider
the set $\boF$ of $t\in\R_+^*$ such, if $\theta_1$ is an extraction, there
exists a sub-extraction $\theta_2$ such that $N_{\theta_2(n)}(Q)$
converge to $N$ at almost every $Q$ of the part of $D$ parametrized by
$\dis]a+\frac{|a|}{t+1},t[$ ($a<0$). Let us prove that
$\boF=\R_+^*$. First, we observe that, if $t_1\in\boF$ and $t_2<t_1$,
$t_2\in\boF$. From Proposition \ref{linofdiv1}, there exists $t>0$
such that $t\in\boF$. Let $t_0=\sup \boF$ and suppose that
$t_0<+\infty$. We consider $P_1$ and $P_2$ the points on $D$
parametrized by $\dis a+\frac{|a|}{1+t_0}$ and $t_0$. We choose $R>0$
such that $D(P_i,R)\subset \Omega$ for $i=1,2$. Let $\theta_1$ be an
extraction. Since $t_0=\sup\boF$, there exist $Q_1\in
D(P_1,\frac{R}{3})\cap D$, $Q_2\in D(P_2,\frac{R}{3})\cap D$ and a
sub-extraction $\theta_2$ such that $N_{\theta_2(n)}$ converges to $N$ at
$Q_1$ and $Q_2$. We have $D(Q_i,\frac{2R}{3})\subset \Omega$, we then
apply Proposition \ref{linofdiv1} to points $Q_1$ and $Q_2$ with
$\alpha=\frac{3}{4}$. We then have a sub-extraction $\theta_3$
such that $N_{\theta_3}$ converges to $N$ at almost every point of
$D(Q_i,\frac{R}{2})\cap D$ for $i=1,2$; this proves that $t_0$ is not
$\sup\boF$, because $N_{\theta_3}$ converges to $N$ at almost every
point of the part of $D$ parametrized by an open interval that
contains the segment $\dis[a+\frac{|a|}{t_0+1},t_0]$.

By a standard diagonal process, we can then construct an extraction
$\theta$ such that $N_{\theta(n)}$ converges to $N$ at almost every point of
$D$. Let $Q$ be in $D$ and we consider $N'$ a cluster point of the
sequence $N_{\theta(n)}(Q)$, if the third coordinate of $N'$ is
negative then there exists a sub-extraction $\theta'$ such that
$W_{\theta'(n)}(Q)$ is bounded but this is impossible since, by Lemma
\ref{bound}, $W_{\theta'(n)}$ would be bounded in a neighborhood of $Q$
and $W_{\theta(n)}$ diverges at almost every point of $D$. Thus the
third coordinate of $N'$ is $0$; if $N'\neq N$, applying what we have
already proved, it appears a second geodesic $D'$ passing by $Q$ normal
to $N'$ and an extraction $\theta'$ such that $N_{\theta'(n)}$
converges to $N'$ at almost every point of $D'$. We parametrized $D$
and $D'$ by arc-length in using the orientation given by the direct
normal to $N$ and $N'$, we choose $Q$ as origin point. Let
$\epsilon>0$; we note $A$ the point on $D$ of coordinate $-\epsilon$
and $B$ the point on $D'$ of coordinate $\epsilon$. For $\epsilon$
small enough, the triangle $AQB$ is in $\Omega$ and then $\dis
\int_{AQB}\dd\Psi_{u_{\theta'(n)}}=0$. We let $n$ tends to $+\infty$
and then obtain $|AC|+|BC|\le |AB|$ which contradicts the triangle
inequality. We then have proved that $N_{\theta(n)}$ converges to $N$
at every point of $D$. We then have proved that for every extraction $\theta$
we can construct a sub-extraction $\theta'$ such that $N_{\theta'(n)}$
converges to $N$ at every point of $D$

To finish the proof, we take a point $Q$ in $D$ and suppose that
$N_{\alpha(n)}(Q)$ converge to $N'$ with $\alpha$ an
extraction. Since $N_{\alpha(n)}(P)\longrightarrow N$, there exists a
sub-extraction $\alpha'$ such that $N_{\alpha'(n)}$ converges to $N$
at every point of $D$, in particular at $Q$, then $N=N'$.
\end{proof}

\begin{rem}{aa}
We then understand the structure of the complementary of $\boB(u_n)$,
it is a set of geodesics of $\Omega$; one of these geodesics will be
called a \emph{line of divergence}. Then when we have a sequence of solutions
of \eqref{MSE}, the problem of the convergence of the sequence is
linked to the understanding of: which lines of divergence are possible?
The answer is, in general, given by the behaviour at the boundary. The
existence of such lines will permit us to use arguments that are
similar to the ones used by H.~Jenkins and J.~Serrin in \cite{JS}. 
\end{rem}

The behaviour of the normal along a line of divergence says us that the
limit of $\dis\int_T\dd\Psi_{u_n}$, where $T$ is a segment of a line of
divergence with the orientation given by the limit normal, is $|T|$. In the
following, we shall draw this limit normal on the figures to explain our
arguments.

\section{The Plateau problem at infinity {\label{ppi}}}

In this section, we explain the problem studied by C.~Cos\'\i n and
A.~Ros in \cite{CR} and give the main results of their paper with some
elements of proofs.

Let M be a complete minimal surface with finite total curvature in
$\R^3$; we know that $M$ is isometric to a compact Riemann surface
minus a finite number of points (we can refer to \cite{Os}). $M$ then has a
finite number of annular ends; when these ends are embedded they are
asymptotic either 
to a half-catenoid or to a plane. A properly immersed minimal surface
with $r$ embedded ends will be called  a \emph{$r$-noid}. We can
associate to each end a vector which caracterizes the direction and the 
growth of the asymptotic half-catenoid (when the end is asymptotic to
a plane this vector is zero); this vector is called the \emph{flux} of
the end (for a precise definition of the flux see \cite{HK}). If
$v_1,\dots,v_r$ are the fluxes at each end, we have the following balancing
condition: 
\begin{equation}{\label{balcond}}
v_1+\cdots+v_r=0
\end{equation}
This condition tells us that the total flux of the system vanishes. If 
$v_1,\dots,v_r$ are vectors in $\R^3$ such that \eqref{balcond} is
verified and $g$ is a non-negative integer, the Plateau problem at
infinity for these data is to find an $r$-noid of genus $g$ which has
$v_1,\dots,v_r$ as fluxes at its ends.

Let $\psi:M\longrightarrow\R^3$  be an $r$-noid. $M$ is conformally
equivalent to a compact surface $\overline{M}$ minus a finite number of
points $p_1,\dots,p_r$. We will say that $M$ is
\emph{Alexandrov-embedded} if $\overline{M}$ bounds a compact $3$-manifold
$\overline{\Omega}$ and the immersion $\psi$ extends to a proper local
diffeomorphism $f:
\overline{\Omega}\backslash\{p_1,\dots,p_r\}\longrightarrow \R^3$. An
Alexandrov-embedded surface has a canonical orientation; we choose the 
Gauss map to be the outward pointing normal. An Alexandrov-embedded
$r$-noid can not have a planar end (see \cite{CR}). We call
$\mathcal{M}_r$ the space of Alexandrov-embedded $r$-noids of genus
$0$ and $r$ horizontal catenoidal ends. We identify two elements  in
$\mathcal{M}_r$ which differ by a translation. In \cite{CR},
C.~Cos\'\i n and A. Ros give a nice description of the space
$\mathcal{M}_r$.

Let $\psi: M\longrightarrow\R^3$ be a nonflat immersion of a connected 
orientable surface $M$ and $\Pi$ be a plane in $\R^3$, normalized to
be $\{x_3=0\}$. We note by $S$ the Euclidiean symmetry with respect to 
$\Pi$ and consider the subsets:

\begin{gather*}
M^+=\{p\in M | x_3(p)>0\}\\
M^-=\{p\in M | x_3(p)<0\}\\
M^0=\{p\in M | x_3(p)=0\}
\end{gather*}
With these notation we have:
\begin{defn}\label{stsym}
We shall say that $M$ is strongly symmetric with respect to $\Pi$ if
\begin{itemize}
\item There exists an isometric involution $s:M\longrightarrow M$
  such that $\psi\circ s=S\circ \psi$.
\item $\{p\in M|s(p)=p\}=M^0$.
\item The third coordinate $N_3$ of the Gauss map of $M$ takes
  positive (resp. negative) values on $M^+$ (resp. $M^-$).
\end{itemize}
\end{defn}

In \cite{CR}, C. Cos\'\i n and A. Ros prove

\begin{prop}
Let $M$ be an $r$-noid with horizontal ends. Then $M$ is strongly
symmetric with respect to an horizontal plane if and only if $M$ is
Alexandrov-embedded.
\end{prop}

We then use the notion of strong symmetry to study $\mathcal{M}_r$; in 
the following, we always suppose that the plane of strong symmetry is
the plane $\{x_3=0\}$. If $M\in \mathcal{M}_r$, $s$ extends to
$\overline{M}=\overline{\C}$, the involution $s$ is $z\mapsto
\dis\frac{1}{\overline{z}}$ and the points $p_1,\dots,p_r$ are in $\{z\in\C |\ 
|z|=1\}$. We then have an order on $\{p_1,\dots,p_r\}$, let us suppose that
$p_1,\dots,p_r$ are put in this order. Let $v_1,\dots,v_r$ be vectors
in $\R^3$ such that $2v_i$ is the flux at the end $p_i$. We have
$v_1+\cdots+v_r=0$, so if we draw the vectors consecutively in the
plane, we get a piece-wise linear closed curve: a polygon. We note
$F(M)$ this polygon.

We say that a polygon $V$ bounds an immersed polygonal disk  if there
exists a compact multi-domain $(\mathcal{P},\phi)$ such that
$\partial\mathcal{P}$ is only composed of edges and
$\phi(\partial\mathcal{P})=V$.

Then the most important result in \cite{CR} is

\begin{thm}
Let $v_1,\dots,v_r$ be horizontal vectors such that
$v_1+\cdots+v_r=0$ and $V$ the associated polygon, then there exists
$M\in \mathcal{M}_r$ such that $F(M)=V$ if, and only if, $V$ bounds an 
  immersed polygonal disk
\end{thm}

Besides, we have as much $M\in\mathcal{M} _r$ such that $F(M)=V$ as
immersed polygonal disks bounded by $V$. Let $V$ be a polygon and
$(\mathcal{P},\phi)$ a compact multi-domain such that
$\phi(\mathcal{P})$ is an immersed polygonal disk bounded by $V$. Let
$P_1,\dots,P_r$ be the vertices of $\boP$ which are identified with the
ones of $V$; we put $P_1=P_{r+1}$. Let $i\in\{1,\dots,r\}$, we can
glue to $\boP$ along $[P_i,P_{i+1}]$ a half-strip $S_i$ isometric
to $[P_i,P_{i+1}]\times\R_+$. We get a multi-domain which we call
$\Ome(\boP)$; the boundary of $\Ome(\boP)$ is composed of $2r$ half
straight-lines, we call $L_i^-$ (resp. $L_i^+$) the half line in the
boundary which has $P_i$ as end point and is in $S_{i-1}$
(resp. $S_i$). 

Let $M$ be in $\boM_r$, we consider $(M^+)^*$ the conjugate surface to
$M^+$ for the outward pointing normal. In \cite{CR}, the authors prove
that it exists $(\boP,\phi)$ a 
multi-domain bounded by $F(M)$ such that $(M^+)^*$ is a graph over the 
multi-domain $\Ome(\boP)$; the normal to the graph is the upward
pointing normal by Definition \ref{stsym}. If $u$ is the function on
$\Ome(\boP)$ that  
gives $(M^+)^*$, they prove that $u$ tends to $+\infty$
(resp. $-\infty$) on $L_i^+$ (resp. $L_i^-$). C. Cos\'\i n and A. Ros
use these arguments to prove that if the Plateau problem at infinity
has a solution the flux polygon $F(M)$ bounds an immersed polygonal
disk. For the other implication, they prove that the map $F: M\mapsto
F(M)$ is a covering map to conclude, they use a compactness argument and
prove that the space $\boM_r$ has a smooth structure.

In the next section, we shall solve on $\Ome(\boP)$ the Dirichlet
problem for the boundary data $+\infty$ on $L_i^+$ and $-\infty$ on
$L_i^-$. We shall then take the conjugate of the graph of the solution
for the downward pointing normal and so build the solution to the
Plateau problem at infinity. The change of orientation makes that we
get the surface we want.

\section{The construction of a solution of the Plateau problem at infinity}

\input{pltinf5.tex}

\appendix

\section{The Carath\'eodory's Theorem{\label{appendice}}}

In this section, we give some explanations on an argument of the proof of
Theorem \ref{regbound}. The problem is: when we have a biholomorphic map
between two open sets of $\C$, can we extend it to the boundary?

We consider $U$ an open set included in $\C$ and $P$ a point of $\partial
U$. We say that $P$ has the property of Sch\"onflies if, for all radii $R$,
there exists a radius $r=r(R)$ such that for all two points in $U\cap D(P,r)$
there exists a path in $U\cap D(P,R)$ joining these two points.

We then have the following theorem that we use in our proof.
\begin{thm}[Carath\'eodory]
Let $U$ be a simply connected open set in $\C$ and $V$ an open set of the
boundary of $U$. We consider $f:U\longrightarrow D=\{z\in\C|\ |z|<1\}$ a
biholomorphic map. We suppose that every point of $V$ has the property of
Sch\"onflies , then $f$ extends to an homeomorphism from $U\cup V$ into $D\cup
C$ where $C\subset\partial D$.
\end{thm}

A proof of this theorem can be found in \cite{He}. In our proof, we have to
verify the property of Sch\"onflies at the points of a part of the
boundary. We know that this part of the boundary is embedded in $\C$ so we
can build neighborhoods of every 
point of the boundary in using $\epsilon$-tubular neighborhood of the
boundary. These neighborhoods prove that we have the property of
Sch\"onflies.

%\begin{figure}[h]
%\begin{center}
%\resizebox{0.5\linewidth}{!}{\input{debase.pstex_t}}
%\end{center}
%\end{figure}

\end{document}

%% file: pltinf0.tex
One classical way to construct minimal surfaces in $\R^3$ is to see
them as the graph of a function $u$ over a domain
$\Ome\subset\R^2$ (see for example the paper of H.~Karcher
\cite{Ka}). The graph of a function $u$ is a minimal surface 
if $u$ satisfies the elliptic partial differential equation called the minimal
surfaces equation:

\begin{equation}{\label{MSE}}
\Div\left(\frac{\nabla u}{\sqrt{1+|\nabla u|^2}}\right)=0
\tag{MSE}
\end{equation}

The problem which is associated to this point of view is the Dirichlet
problem for the equation \eqref{MSE}: for a domain $\Omega$ and a
function $f$ on 
$\partial\Ome$, this problem consits in finding a continuous function
$u$ on $\overline\Ome$ which is a solution of the minimal surfaces
equation in $\Ome$ and such that $u=f$ on the boundary of $\Ome$. One 
of the most  general answers to the Dirichlet problem for bounded
domain has been given by H. Jenkins and J. Serrin in \cite{JS}. They
give a nice condition on the domain to solve for any function $f$;
moreover, their
result allows us to give infinite value for the
boundary data $f$.  For unbounded domain, the Dirichlet problem is
still an open problem. We know that, in the general case, we lose the
uniqueness of solution. In this paper, using a new approach, we
develop some tools for the study of this problem. 

An other interesting and still open problem concerning minimal
surfaces is the Plateau problem at 
infinity which is the following: finding a minimal surface for 
a given
asymptotic behaviour. More precisely, we know that, if a complete
minimal surface has finite 
total curvature and embedded ends, each end of this minimal surface
is asymptotic to a plane or to a half-catenoid; besides, we can
associate to each end a vector in $\R^3$, this vector is called the
flux vector of the end. These vectors satisfy the following condition:
the sum of the flux vectors over all ends is zero. So the problem is: given a
finite number of vectors such that their sum is zero, can we find a
minimal surface which has these vectors as flux vectors? Our answer
comes from the following idea: seeing a solution of the Plateau problem
at infinity as the conjugate surface of a solution of the Dirichlet
problem on an unbounded domain. 

In \cite{CR}, C. Cos\'\i n and A. Ros give a description of the space of
solutions of the Plateau problem at infinity with an asymptotic
behaviour which 
is symmetric with respect to an horizontal plane (\emph{i.e.} all the
flux vectors are horizontal). They also restrict themselves to the
case of Alexandrov embedded minimal surfaces; this condition implies
that no flux vector is zero and that there is a natural order on the
ends of the surface. Since the flux vectors
are horizontal 
and their sum is zero, these vectors draw a polygon in
$\R^2$. C.~Cos\'\i n and A.~Ros give a necessary and
sufficient condition on this polygon to have a solution. See section
\ref{ppi}, for more explanations about their work. 

In this paper, we give a more constructive proof of the result of C. Cos\'\i n
and A. Ros. Our method is based on the Dirichlet problem on an
unbounded ``domain'' $\Ome$. When 
the polygon given by the flux vectors is convex, $\Omega$ can be
defined as the polygonal domain bounded by the flux polygon to which we glue a
half-strip on each edge. We note $L_i^+$ and $L_{i+1}^-$ the two sides of each
half-strip $S_i$, alternating the sign $+$ and $-$ such that each
vertex of the polygon is common to some $L_i^-$ and $L_i^+$. When the
flux polygon is non-convex and satisfies the condition 
of C.~Cos\'\i n and A.~Ros, we need to use the concept of mutli-domain
for defining $\Ome$ (see Definition \ref{multdom} for this concept).

Our main result for the Dirichlet problem for this kind of domain
$\Ome$ is then (see Theorem \ref{moi}): 

\bigskip

\emph{There exists a solution $u$ of the minimal surfaces   equation
  on $\Ome$ such that $u$ tends to $+\infty$ on $L_i^+$ and $-\infty$
  on $L_i^-$. Besides, the solution is unique up to an additive
  constant} 
 
\bigskip

\noindent
The function $u$ in this result is build as the limit of solutions of
the Dirichlet problem on bounded domain. We describe
the possible divergences that can occur for a sequence of
solutions of 
\eqref{MSE}. In fact, we prove that if the sequence diverges at a
point, it must diverge along a line passing by this point. This result
is a generalization of the results that H.~Jenkins and J.~Serrin use
in \cite{JS}. Our result allows us to do the same discussion that
H.~Jenkins and J.~Serrin made in the particular case of monotone
sequences of 
solutions of \eqref{MSE}; this is our main tools to prove the
existence part of Theorem \ref{moi}.

The solution to the Plateau problem at infinity is then the conjugate
surface to the graph of $u$. In order to know the geometry of the conjugate
surface along its boundary, we need to understand the behaviour of the
graph in the 
neighborhood of the vertices of $\Ome$ which are the vertices of the
polygon. Some results are known for such problem concerning the
Dirichlet problem in the convex case. For example, consider
$f$ a data on the boundary of a domain $\Ome$, we suppose that $f$ has
a finite discontinuity 
at a point $P$ where the boundary is convex (\emph{i.e.} we suppose
that $f(Q)$ has a limit if we tend to $P$ by the right hand side or by
the left hand side and that the difference of these two limits is
finite), then we know that the graph of a solution $u$ over $\Ome$ of
the 
Dirichlet problem with $f$ as boundary value, has a vertical segment
over $P$ in its boundary, it was proved in \cite{Ni2}. In our case, we
can prove that the boundary of 
the graph is the vertical straight line passing by the vertex
; although the domain is locally an angular sector that not need to be
convex and 
the boundary data takes the values $+\infty$ on one side of the sector
and $-\infty$ on the other side.

\bigskip

The paper is organized as follows; in the first section, we define multi-domains and extend the result of
H.~Jenkins and J.~Serrin to bounded multi-domains. The multi-domains
are necessary to express the condition of C.~Cos\'\i n and
A.~Ros. This result will be our first tool in the proof of  our main
theorem.

The second section is devoted to the proof of our result concerning
the boundary behaviour of solutions of the Dirichlet problem.

In section \ref{seqmse}, we study the sequences of solutions of
\eqref{MSE} and define the lines of divergence.

In section \ref{ppi}, we explain the result of C.~Cos\'\i n and A.~Ros, and
recall some elements of their proof. In the last section, we give the proof of
our main result. We then use it to give a new proof of the result of C.~Cos\'
\i n and A.~Ros. 
 
Let us fix some notations. In the following, when $u$ is a function on a
domain of $\R^2$ we shall note $W=\sqrt{1+|\nabla u|^2}$. We shall also use
the classical following notations for partial derivatives: $p=\der{u}{x}$,
$q=\der{u}{y}$, $r=\dis\frac{\partial^2u}{\partial x^2}$,
$s=\dis\frac{\partial^2 u}{\partial x\partial y}$ and $t=\dis \frac{\partial^2
  u}{\partial y^2}$. Besides, for the graph of $u$, we shall always chose the
downward pointing normal to give an orientation to the graph.

\bigskip

%\textbf{Acknowledgements.} The author would like to thank Pascal Collin
%for many helpful discussions. 

%% file: pltinf2.tex
The aim of this section is to understand what geometrically happens at
a vertex of a multi-domain where two edges $A_j$ and $B_j$ converge.

For $\beta_1<\beta_2$ and $R>0$, we consider:
$$\Ome_{\beta_1}^{\beta_2}(R)= \left\{(r,\theta)|\ 0\le r\le R,\ \beta_1\le
\theta\le \beta_2\right\}$$
with the metric $\dd s^2=\dd r^2+r^2\dd \theta^2$ (we identify all the points
$(0,\theta)$ and this point will be called the vertex of
$\Ome_{\beta_1}^{\beta_2}(R)$). We define also on
$\Ome_{\beta_1}^{\beta_2}(R)$ the map $\phi:(r,\theta) \mapsto
(r\cos\theta,r\sin\theta)$. Then $(\Ome_{\beta_1}^{\beta_2}(R),\phi)$ is a
multi-domain, it is a description of a neighborhood of a vertex where
two edges converge. We call $L(\beta)$ the set of points in
$\Ome_{\beta_1}^{\beta_2}(R)$ such that $\theta=\beta$. We are interested in
the geometrical ``configuration'' of the graph of a solution $u$ of
\eqref{MSE} such that $u$ tends to $-\infty$ on $L(\beta_2)$ and $+\infty$ on
$L(\beta_1)$; such a solution $u$ will be called a solution of the problem
$\boP$. 

The first thing we have to do to understand a solution $u$ of the
problem 
$\boP$ is being able to bound the function $u$ on each radius
$L(\beta)$. Our arguments are based on the comparison with the Scherk
surface. 

Let us consider $ABC$ an isosceles triangle ($|AB|=|AC|=R$), we consider
the solution $w$ of the Dirichlet problem on $ABC$ such that $w=0$ on
$[A,B]$ and $[A,C]$ and tends to $+\infty$ on $[B,C]$; this function
exists by Theorem \ref{thmjs}. When $ABC$ is
rectangle $w$ is the Scherk surface, after dilatation, $w$ is given by:

\begin{equation}
w(x,y)=h(x,y)=-\ln\cos x+\ln\cos\ y
\end{equation}
In the general case, the solution $w$ will be called a 
\emph{pseudo Scherk surface}. 

We shall use the Scherk surface to control solutions of the
problem $\boP$. We first consider the case where $ABC$ is rectangle. In fact,
a neighborhood of $B$ in $ABC$ can be
isometrically parametrized by $\Ome_{-\frac{\pi}{4}}^0(R)$ and $h$ is a
solution of \eqref{MSE} on $\Ome_{-\frac{\pi}{4}}^0(R)$ such that $h=0$ on
$L(0)$, $+\infty$ on $L(-\frac{\pi}{4})$ and some positive function on the
third part of the boundary. Since we have an expression for $h$ we can 
see that $h$ is uniformly bounded on $\Ome_{\alpha}^0(R)
\subset\Ome_{-\frac{\pi}{4}}^0(R)$ for every
$-\frac{\pi}{4}<\alpha<0$.

We do not suppose now that $ABC$ is rectangle; but we suppose that the 
angle at the vertex $A$ is greater than $\frac{\pi}{2}$. In this case we can
choose a point $A'$ such that $A'BC$ is isosceles and rectangle and $A'BC$
contains $ABC$. We consider in
$ABC$ the pseudo Scherk surface $w$ and $h$ the Scherk surface on $A'BC$;
since $h$ is positive in $A'BC$, we have $h>w$. As above, a neighborhood of $B$
in $ABC$ can be isometrically parametrized by $\Ome_\beta^0(R)$ with
$\beta<0$ and $w$ can be seen as the solution of \eqref{MSE} on
$\Ome_\beta^0(R)$ such that $w=0$ on $L(0)$, $+\infty$ on $L(\beta)$
and some positive function on the third part of the boundary. Since $w<h$, $w$
is uniformly bounded on $\Ome_\alpha^0(R)$ for every $\beta<\alpha<0$.

By our expression for $h$, there exists $m\in\R$ such that $h\le m$ on $[A,B]$
and $[A,C]$. This proves that $h-m\le w$ in $ABC$. Then in our parametrization
of a neighborhood of $B$, for every $M\in\R$ there exist $\alpha$ such that
$w\ge M$ in $\Ome_{\beta}^\alpha(R)$. 

\begin{lem}{\label{psherck}}
Let $\beta_1<\beta_2$ and $R>0$. We consider a solution $u$ of the problem
$\boP$ on $\Ome_{\beta_1}^{\beta_2}(R)$. Then for every
$\beta_1<\alpha<\beta_2$, there exist $M$ and $M'$ in $\R$ such that $u\le M$
in $\Ome_\alpha^{\beta_2}(\frac{R}{4})$ and $u\ge M'$ on
$\Omega_{\beta_1}^\alpha(\frac{R}{4})$. For every $M\in\R$, there exist
$\alpha$ and $\alpha'$ in $]\beta_1,\beta_2[$ such that $u\ge M$ in
$\Ome_{\beta_1}^\alpha(\frac{R}{4})$ and $u\le M$ in $\Ome_{\alpha'}^{\beta_2}
(\frac{R}{4})$.
\end{lem}

\begin{proof}
Let us consider $\alpha_1<\alpha_2$ and $R'>0$. We consider $v$ the solution of
the problem $\boP$ on  $\Ome_{\alpha_1}^{\alpha_2}(R')$ such that $v=0$ on the
third part of the boundary; $v$ exists because the hypotheses of
Theorem \ref{thmjs} are fulfilled. The isometry of
$\Ome_{\alpha_1}^{\alpha_2}(R')\times \R$
defined by $(r,\theta,z)\mapsto(r,\alpha_1+\alpha_2-\theta,-z)$ does not
change the
boundary data so $v$ is invariant by this isometry because of the uniqueness
of such a solution. This proves that $v=0$ on
$L(\frac{\alpha_1+\alpha_2}{2})$. Then, by maximum principle, we have $v>0$
between $L(\alpha_1)$ and $L(\frac{\alpha_1+\alpha_2}{2})$ and $v<0$ between
$L(\frac{\alpha_1+\alpha_2}{2})$ and $L(\alpha_2)$.

Let us consider $\alpha_1<\alpha<\alpha_2$. Let us prove that there
exists a constant $M$ such that $v\le M$ in
$\Ome_{\alpha}^{\alpha_2}(\frac{R'}{2})$. If $\alpha\ge
\frac{\alpha_1+\alpha_2}{2}$, $M=0$ works. We note
$\widetilde{\alpha}=\frac{\alpha_1+\alpha_2}{2}$. We suppose $\alpha\le
\widetilde{\alpha}$ then we take a sufficiently big $n$ such that
$\frac{\widetilde{\alpha} -\alpha}{n}\le \frac{\pi}{4}$ and
$\frac{\widetilde{\alpha} -\alpha}{n}\le \alpha-\alpha_1$. We note $B$ the
vertex of $\Ome_{\alpha_1}^{\alpha_2}(R')$. For
$k\le 2n+1$ we note $\alpha(k)=\dis\widetilde{\alpha}
-k\frac{\widetilde{\alpha}  -\alpha}{2n}$ and for $k\le 2n-1$ we note $A_k$
the points of coordinates 
$(\frac{R'}{2},\alpha(k))$ and $C_k$ the point of second coordinate
$\alpha(k+2)$ such that $A_kBC_k$ is an isosceles triangle at $A_k$ (see
Figure \ref{triangle}). We have $v=0$ on
$[B,A_0]$ and $v$ is bounded on $[A_0,C_0]$, then, by adding a constant, we
can put a pseudo Scherk surface above $v$ over $A_0BC_0$. This proves that $v$
is upper-bounded in $\Ome_{\alpha(1)}^{\alpha_2}(\frac{R'}{2})$. Since $v$ is
upper-bounded on $[B,A_1]$ and $[A_1,C_1]$, we can put a pseudo Scherk surface
above $v$ over $A_1BC_1$ then $v$ is bounded on
$\Ome_{\alpha(2)}^{\alpha_2}(\frac{R'}{2})$. We can do this for every $k$ then
we obtain that $v$ is uniformly upper-bounded on
$\Ome_{\alpha}^{\alpha_2}(\frac{R'}{2})$.

\begin{figure}[h]
\begin{center}
\resizebox{0.8\linewidth}{!}{\input{figure3.pstex_t}}
\caption{\label{triangle}}
\end{center}
\end{figure}

With the same method, we can prove that there exists $M'$ such that $v\ge M'$
on $\Ome_{\alpha_1}^{\alpha}(\frac{R'}{2})$.

Let us now consider our original problem. We have $u$ and $\alpha$ and we want
to prove the existence of $M$. We consider $\beta_1<\alpha'<\alpha$, since
$u$ tends to $-\infty$ along $L(\beta_2)$, there exists $m$ such that $u\le
m$ at all the points of coordinates $(\frac{R}{2},\theta)$ with
$\alpha'\le\theta<\beta_2$. We consider on
$\Ome_{\alpha'}^{\beta_2}(\frac{R}{2})$ the solution $v$ that we have studied
above, by maximum principle, we have $u\le v+m$ on
$\Ome_{\alpha'}^{\beta_2}(\frac{R}{2})$. We then have the existence of $M$
because of the result on $v$. We construct $M'$ in the same way. 

Let us now consider $u$ a solution of $\boP$ and $M\in\R$. We consider
$\beta_1<\beta<\beta_2$ such that $\beta-\beta_1\le \frac{\pi}{4}$, we
consider the point $A$ of coordinates $(\frac{R}{2},\beta)$ and $C$ the point
on $L(\beta_1)$ such that $ABC$ is a isosceles triangle (where $B$ is the
vertex of $\Ome_{\beta_1}^{\beta_2}$). By what we have just proved, $u$ is
lower-bounded on $[B,A]$ and $[A,C]$ then we can put a pseudo Scherk
surface under $u$. The existence of $\alpha$ is due to the last remark
that we made about pseudo Scherk surfaces
\end{proof}

Using this result, we can prove the following geometrical result.

\begin{thm}{\label{regbound}}
Let $(\Ome, \phi)$ be a multi-domain and $P$ a vertex of $\Ome$ such that two
edges $L_1$ and $L_2$ have $P$ as end point ($L_1$ and $L_2$ are enumerated
with respect to  the orientation). Let $u$ be a solution of \eqref{MSE} on
$\Ome$ such that $u$ tends to $-\infty$ on $L_1$ and $+\infty$ on $L_2$. We
consider $\Psi_u$ the conjugate function to $u$ normalized such that
$\Psi_u(P)=0$. Then, if $\Psi_u$ is non-negative in a neighborhood of
$P$, the 
vertical straight line passing through $\phi(P)$ is the boundary of the graph
of $u$ above a neighborhood of $P$. 
\end{thm}

First, we remark that, if $Q$ is a point on $L_1$ or $L_2$, then
$\Psi_u(Q)=|PQ|\ge 0$ by Lemma \ref{calc}. This proves that, if the
angle at $P$ is strictly less 
than $\pi$, the hypothesis on $\Psi_u$ is always verified; so we have the
result for a convex corner.

\begin{proof}
By a translation and a rotation, we can isometrically parametrized a
neighborhood of $P$ by $\Ome_{\beta}^0(R)$ for $\beta<0$ and $R$ small
enough. Then $u$ can be seen as a solution of the problem $\boP$. We suppose
that $\Psi_u\ge 0$ in $\Ome_\beta^0(R)$.

\begin{center}
\underline{First part}
\end{center}

First, we prove that there exists $M_1\in\R$ such that
$\phi(P)\times]-\infty,M_1[$ is a part of the boundary of the graph. We take
$-\frac{\pi}{2}<\alpha<0$, we suppose that $\alpha>\frac{\beta}{2}$. Then
$\Ome_{\alpha}^0(R)\subset\Ome_\beta^0(R)$ can be parametrized by euclidean
parameters $(x,y)$, in fact $\Ome_\alpha^0(R)$ is embedded in $\R^2$. The idea
is to see the part of the graph which is over $\Ome_\alpha^0(R)$ as a
graph over the vertical plane given by the equation $y=0$. Let
$R'<\frac{R}{2}$, then for all 
$Q\in \Ome_{\alpha}^0(R')$ the nearest point from $Q$ on $\partial\Ome$ is on
$L_1$. If we take $R'$ small enough and $\alpha$ such that
$\tan\alpha>-\frac{1}{8}$, then every point of $\Ome_\alpha^0(R')$ verifies the
hypothesis of Lemma 1 in \cite{JS}. This lemma implies that, at every point of
$\Ome_\alpha^0(R')$, $q=\der{u}{y}<0$. Using our euclidean parameters,
we note, for $(x,y)\in\Ome_\alpha^0(R')$, $\Theta(x,y)=(x,u(x,y))$. We have:
$$\dd\Theta=
\begin{pmatrix}
1&0\\
\der{u}{x}&\der{u}{y}
\end{pmatrix}
$$
Since $q<0$, this proves that $\Theta$ is a local
diffeomorphism. Since $u$ strictly 
decreases when $y$ increases, $\Theta$ is injective. By Lemma
\ref{psherck}, we know that there exists
$K\in\R$ such that $u\ge K$ on $L(\alpha)$, we put $x_1=R'\cos\alpha$. We then
have $]0,x_1[\times]-\infty,K[\subset\Theta\left(\Ome_\alpha^0(R')\right)$. We
note $\chi=\Theta^{-1}$ on $]0,x_1[\times]-\infty,K[$; we then have
$y=\chi_{(2)}(x,z)$ on the graph of $u$ (we note $\chi_{(2)}$ the second
coordinate function of $\chi$) then $\chi_{(2)}$ verifies
\eqref{MSE}. When $x\longrightarrow 0$, we have $y=\chi_{(2)}(x,z)
\longrightarrow 0$, it is due to the shape of $\Ome_\alpha^0(R')$. From Lemma
\ref{psherck}, there exist $\alpha<\alpha'<0$ and $r$ such that
$\Ome_{\alpha'}^0(r)\subset \im\chi$. By results of boundary regularity, 
$\chi_{(2)}$ is regular at the boundary, actually we can extend
$\chi_{(2)}$ by making a reflection with respect to the axis
$\{x=0,y=0\}$. We now show that a part of this axis is a part of the
boundary of the whole graph.
By lemma \ref{psherck}, there exists 
$M'$ such that $u\ge M'$ in $\Ome_{\beta}^{\alpha'}(r)$. We note $M_1=M'-1<K$;
then if a sequence of points of the graph of $u$ over $\Ome_\beta^0(r)$ tends
to a point of $\phi(P)\times]-\infty,M_1[$, we have $(x,y)$ in $\im\chi$ after
a certain rank. Then the graph of $u$ over $\im\chi$ is a neighborhood of
$\phi(P)\times]-\infty,M_1[$; as $\chi_{(2)}$ is  regular through the
boundary, $\phi(P)\times]-\infty,M_1[$ is a part of the boundary.

With the same arguments, we can show that there exists $M_2$ such
that $\phi(P)\times 
]M_2,+\infty[$ is a part of the boundary.

\begin{center}
\underline{Second part}
\end{center}

The first part proves that outside a compact the graph of $u$ has a
good behaviour above the point $\phi(P)$. Now, we prove that we can
extend, by reflection, this compact part through the verical straight
line passing by $\phi(P)$. 

From what we have just done, there exist $\beta<\alpha_2<\alpha_1<0$ such that
the graph above $\Ome_\beta^{\alpha_2}(R)$ and $\Omega_{\alpha_1}^0(R)$ is
regular above $P$. We choose $M_1$ and $M_2$ as in the first part such that
$(\phi(P), M_1)\in \partial\Gra(u|_{\Ome_{\alpha_1}^0(R)})$ and $(\phi(P),
M_2)\in \partial\Gra(u|_{\Ome_\beta^{\alpha_2}(R)})$. We shall construct a
curve $\Gamma$ as follow: we start from the point $A_1=(0,0,M_1)$
($(0,0)=\phi(P)$), we go down vertically to the point
$A_2=(0,0,M_1-1)$, then 
we go to some point $A_3=(\epsilon\cos\theta, \epsilon\sin\theta, M_1-1)$ in
following the level curve $\{u=M_1-1\}$ (we suppose $\epsilon$ small and
$\theta>\alpha_1$), we then follow  the curve 
$$t\mapsto (\epsilon\cos t ,\epsilon \sin t, u(\epsilon\cos t ,\epsilon \sin
t)),$$
we let $t$ decreases to some $\theta'<\alpha_2$ such that $u(\epsilon\cos
\theta' ,\epsilon \sin\theta')=M_2+1$ (we note $A_4$ the end point), following
the level curve $\{u=M_2+1\}$ we go to the point $A_5=(0,0,M_2+1)$ and finally
we go down to the point $A_6=(0,0,M_2)$. We can smooth $\Gamma$ at the points
$A_2,A_3,A_4$ and $A_5$ such that $A_2$ and $A_5$ are always in the smooth
$\Gamma$ and the new $\Gamma$ is embedded in the graph of $u$. The vertical
projection of $\Gamma$ on $\Ome_\beta^0(R)$ bounds a domain
$\widetilde{\Ome}$. We note $\Sigma$ the graph of $u$ above
$\widetilde{\Ome}$. Because of our choice of $\Gamma$, $\Sigma$ extends in a
minimal surface $\Sigma'$ through $\Gamma$ (The only problem is
through $[A_1,A_2]$ and $[A_5,A_6]$, but the first part says us that
we can extend $\Sigma$ through these two segments by symmetry). Because
$\Sigma$ is 
a graph, $\Sigma$ is simply connected and its boundary is not empty; the same
is 
true for $\Sigma'$. This remark says us that we have conformal parametrization
$h_1:D\longrightarrow \Sigma'$ and $h_2:D\longrightarrow\Sigma$ ($D$ is the
unit disk). We put $\widetilde{D}=h_1^{-1}(\Sigma)$ and
$\widetilde{h}:\widetilde{D}\longrightarrow D$ defined by $\widetilde{h}=
h_2^{-1} \circ h_1$; $\widetilde{h}$ is a biholomorphic map. As
$h_1^{-1}(\Gamma)$ is embedded in $D$, the property of Sch\"onflies is
verified at every point ; by the Carath\'eodory's Theorem,
$\widetilde{h}$ extends 
to an homeomorphism of $\widetilde{D} \cup h_1^{-1}(\Gamma)$ into $D\cup V$
where $V$ is part of the boundary of $D$ (for all this argument we refer to
appendix \ref{appendice}). This 
proves that we can extend $h_2$ 
in an homeomorphism of $D\cup V$ into $\Sigma\cup \Gamma$. Let us consider
$f:D\longrightarrow D^-$ ($D^-=\{(x,y)\in D|\ y<0\}$) a biholomorphic map,
then $f$
extends to the boundary. Let us consider the following points on $\Gamma$:
$A_{1.5}=(0,0,M_1-0.5)$ and $A_{5.5}=(0,0,M_2+0.5)$. We note $X=h_2\circ
H\circ f^{-1}$ where $H$ is a Moebius transformation of the unit disk. We note
$B_i=X^{-1}(A_i)$ for every $i$. Then, for a suitable choice of $H$, we can
have the situation described by Figure \ref{reg}.

\begin{figure}[h]
\begin{center}
\resizebox{0.8\linewidth}{!}{\input{figure2.pstex_t}}
\caption{\label{reg}}
\end{center}
\end{figure}

Let us show that $X$ extends to the whole disk. We shall note
$x_1,x_2$ and $x_3$ the three coordinates of $X$, this three functions are
harmonic since $\Sigma$ is minimal. First we observe that $x_1$ and $x_2$ tend
to $0$ when $z\in D^-$ 
tends to $D^0=\{z\in D|\ z\in\R\}$, this is due to the shape of
$\widetilde{\Ome}$. Then, by Schwarz reflection principle,
$x_1$ and $x_2$ extend to $D$ in harmonic functions. Let us consider $x_3^*$
the harmonic conjugate to $x_3$ on $D^-$, we normalized $x_3^*$ by
$x_3^*(B_2)=0$. By our choice of normalization, for every $z\in D^-$, we have
$\Psi_u(X(z))=x_3^*(z)$; this proves that $x_3^*$ tends to $0$ when $z$ tends
to $D^0$, we can extend $x_3^*$ by reflection to $D$. By taking the conjugate
function, we have proved that $x_3$ extends to $D$. We then have
constructed a minimal immersion $X$ on $D$, then $\Sigma$ extends
through $[A_1,A_6]$. This extention is given by the reflection with respect to 
$[A_1,A_6]$.

\begin{center}
\underline{Third part}
\end{center}     

The last thing we have to show is that the minimal immersion $X$ has no branch
point. If it has a branch point then it must be on $D^0$, since, on the other
part, the surface is a graph and then there is no branch point. Let $z_0$ be a
branch point, then $\nabla x_3^*(z_0)=0$. Since $x_3^*$ is harmonic, its local
behaviour is quite similar to the one of $\Re(z-z_0)^p$ with $p\ge 2$ (in
fact, in some holomorphic chart we have $x_3^*(z)=\Re(z-z_0)^p$). This
implies that there exists $z$ in $D^-$ such that $x_3^*(z)<0$, but this
contradicts our hypothesis $\Psi_u\ge 0$.  

We then have proved that there is no branch point so the vertical straight
line passing by $\phi(P)$ is the boundary of the graph.
\end{proof}

\begin{rem}{qjfsd}
We can remark that in the first two parts we do not use the hypothesis on
$\Psi_u$. So in such a situation we can always extend the graph by making a
reflection with respect to the vertical axis. But what we obtain is a minimal
surface with, may be, a finite number of branch points on the vertical axis. 
\end{rem}

\begin{rem}{cc}
We can make an other remark. If we consider a vertex $P$, two edges 
$L_1$ and $L_2$ having $P$ as end point, and $u$ such that $u$ assumes
the data $+\infty$ (or $-\infty$) on $L_1$ and $L_2$, the hypothesis
on $\Psi_u$ did not have any more sense and the angle at the vertex
$P$ must be greater than $\pi$. But we can always apply the
two first parts of the proof. The only problem is that we need a
result similar to Lemma \ref{psherck}; this is given by Theorem 10.3
in \cite{Os}. So we can affirm that on the boundary of the graph of
$u$ we have a half straight line with a finite number of branch
points. Obviously, we must have a branch point at the end point of the 
half straight line.
\end{rem}

%%% Local Variables: 
%%% mode: latex
%%% TeX-master: t
%%% End: 

%% file: figure3.pstex_t
\begin{picture}(0,0)%
\epsfig{file=figure3.pstex}%
\end{picture}%
\setlength{\unitlength}{4144sp}%
\begingroup\makeatletter\ifx\SetFigFont\undefined%
\gdef\SetFigFont#1#2#3#4#5{%
  \reset@font\fontsize{#1}{#2pt}%
  \fontfamily{#3}\fontseries{#4}\fontshape{#5}%
  \selectfont}%
\fi\endgroup%
\begin{picture}(7482,6133)(889,-7838)
\put(4726,-6181){\makebox(0,0)[lb]{\smash{\SetFigFont{12}{14.4}{\rmdefault}{\mddefault}{\updefault}$A_{2n-1}$}}}
\put(901,-6451){\makebox(0,0)[lb]{\smash{\SetFigFont{12}{14.4}{\rmdefault}{\mddefault}{\updefault}$L(\alpha_1)$}}}
\put(4051,-7171){\makebox(0,0)[lb]{\smash{\SetFigFont{12}{14.4}{\rmdefault}{\mddefault}{\updefault}$L(\alpha)$}}}
\put(2521,-7531){\makebox(0,0)[lb]{\smash{\SetFigFont{12}{14.4}{\rmdefault}{\mddefault}{\updefault}$C_{2n-1}$}}}
\put(3376,-1861){\makebox(0,0)[lb]{\smash{\SetFigFont{12}{14.4}{\rmdefault}{\mddefault}{\updefault}$L(\alpha_2)$}}}
\put(6841,-3481){\makebox(0,0)[lb]{\smash{\SetFigFont{12}{14.4}{\rmdefault}{\mddefault}{\updefault}$A_0$}}}
\put(6751,-4066){\makebox(0,0)[lb]{\smash{\SetFigFont{12}{14.4}{\rmdefault}{\mddefault}{\updefault}$A_1$}}}
\put(6571,-4561){\makebox(0,0)[lb]{\smash{\SetFigFont{12}{14.4}{\rmdefault}{\mddefault}{\updefault}$A_2$}}}
\put(8371,-5506){\makebox(0,0)[lb]{\smash{\SetFigFont{12}{14.4}{\rmdefault}{\mddefault}{\updefault}$C_0$}}}
\put(7246,-6811){\makebox(0,0)[lb]{\smash{\SetFigFont{12}{14.4}{\rmdefault}{\mddefault}{\updefault}$C_2$}}}
\put(4771,-3526){\makebox(0,0)[lb]{\smash{\SetFigFont{12}{14.4}{\rmdefault}{\mddefault}{\updefault}$B$}}}
\put(7876,-6181){\makebox(0,0)[lb]{\smash{\SetFigFont{12}{14.4}{\rmdefault}{\mddefault}{\updefault}$C_1$}}}
\end{picture}

%% file: figure2.pstex_t
\begin{picture}(0,0)%
\epsfig{file=figure2.pstex}%
\end{picture}%
\setlength{\unitlength}{4144sp}%
\begingroup\makeatletter\ifx\SetFigFont\undefined%
\gdef\SetFigFont#1#2#3#4#5{%
  \reset@font\fontsize{#1}{#2pt}%
  \fontfamily{#3}\fontseries{#4}\fontshape{#5}%
  \selectfont}%
\fi\endgroup%
\begin{picture}(3915,2400)(1531,-4111)
\put(1846,-1906){\makebox(0,0)[lb]{\smash{\SetFigFont{12}{14.4}{\rmdefault}{\mddefault}{\updefault}$B_{1.5}$}}}
\put(1531,-2626){\makebox(0,0)[lb]{\smash{\SetFigFont{12}{14.4}{\rmdefault}{\mddefault}{\updefault}$B_2$}}}
\put(3781,-4111){\makebox(0,0)[lb]{\smash{\SetFigFont{12}{14.4}{\rmdefault}{\mddefault}{\updefault}$B_4$}}}
\put(5446,-2671){\makebox(0,0)[lb]{\smash{\SetFigFont{12}{14.4}{\rmdefault}{\mddefault}{\updefault}$B_5$}}}
\put(4816,-1951){\makebox(0,0)[lb]{\smash{\SetFigFont{12}{14.4}{\rmdefault}{\mddefault}{\updefault}$B_{5.5}$}}}
\end{picture}

%% file: pltinf5.tex
The first part of this section will be devoted to the proof of our
main result. 

\begin{thm}\label{moi}
Let $V$ be a polygon wich bounds an immersed polygonal disk
$(\boP,\phi)$, we define $\Ome(\boP)$ as in the preceding section. Then
there exists a solution $u$ of \eqref{MSE} on $\Ome(\boP)$
such that $u$ tends to $+\infty$ on $L_i^+$ and $-\infty$ on
$L_i^-$. Besides, the solution is unique up to an additive constant.
\end{thm}

Let us first consider $u$ a solution of \eqref{MSE} on the half-strip
$[0,a]\times\R_+$ such that $u$ tends to 
$-\infty$ on $\{a\}\times\R_+^*$ and $+\infty$ on
$\{0\}\times\R_+^*$. This situation describes the behaviour in the $r$
half-strips $S_i$. Then by Lemma $1$ in \cite{JS} we have:

\begin{gather}
\frac{|q|}{W}(x,y)\ge 1-\frac{a^2}{x^2}\\
\frac{|p|}{W}(x,y)\le \sqrt{2}\frac{a}{x}{\label{bb}}
\end{gather}
when $x\ge 4a$. We consider now the general problem.

We begin in proving the uniqueness part of Theorem \ref{moi}. Let
$u_1$ and $u_2$ be two 
different solutions of the problem (\emph{i.e.} $u_1-u_2$ is
non-constant). As in the proof of Theorem \ref{thmjs}, we can suppose
that $\{u_1>u_2\}$ and $\{u_1<u_2\}$ are non-empty. Let us call
$\Ome_l$ the subset of $\Ome(\boP)$ which is the union of $\boP$ and
the set of points in each $S_i$ that are at a distance less than $l$ from
$[P_i,P_{i+1}]$; we define $\Ome_l^+=\Ome_l\cap\{u_1>u_2\}$. Let us
consider:
\begin{equation}
I=\int_{\partial\Ome_l^+}\dd \tilde{\Psi}
\end{equation}
where $\dd \tilde{\Psi}=\dd \Psi_{u_1}-\dd \Psi_{u_2}$. Since $\dd
\tilde{\Psi}$ is closed, we have $I=0$. $\partial\Ome_l^+$ is composed 
of a part which is included on $\left( \cup_i L_i^+]\right) \bigcup \left
  ( \cup_i L_i^-\right)$ where $\dd \tilde{\Psi}=0$, a part included in the
interior of $\Ome(\boP)$, 
noted $\Gamma_l$, and a part in $I_{i,l}$ which is the part in
$S_i$ parametrized by $[P_i,P_{i+1}]\times\{l\}$. On the part
included in $I_{i,l}$ if $l$ is big enough the integral of $\dd
\tilde{\Psi}$ is less than $2\sqrt{2}\frac{|P_iP_{i+1}|^2}{l}$ by
\eqref{bb}. We then have:

\begin{equation}
0=I\le \int_{\Gamma_l}\dd\tilde{\Psi} +\sum_{i=1}^r 2\sqrt{2}
\frac{|P_iP_{i+1}|^2}{l}
\end{equation}      

By Lemma $2$ in \cite{CK}, $\dis\int_{\Gamma_l}\dd\tilde{\Psi}$ is
negative and decreases as $l$ increases. Because $\dis
\sum_{i=1}^r \sqrt{2}\frac{|P_iP_{i+1}|^2}{l}
\xrightarrow[l\rightarrow +\infty]{} 0$, we get a contradiction. This
proves that, if $u_1$ and $u_2$ are two solutions of our Dirichlet
problem, there exists $c\in\R$ such that $u_1=u_2+c$.

We now prove the existence of the solution. We fix a point $P_0$ in
$\boP$. Let us consider in $S_i$ the point $Q_i^k$ which is the middle
point of $I_{i,k}$, we then define $\Ome_k$ to be  the compact
subdomain of $\Ome(\boP)$ bounded by the segments $[P_i,Q_i^k]$ and
$[Q_i^k,P_{i+1}]$. Let $\boG_i^k$ be  the set of the points $Q$ in
$\Ome_k$ such that $d(Q,Q_i^k)<d(P_i,Q_i^k)$; if $k$ is big enough the
sets $\boG_i^k$ are disjoint, this proves that the conditions of Theorem
\ref{thmjs} are fulfilled for big $k$. Then by Theorem~\ref{thmjs}, we can
build a function $u_k$ on $\Ome_k$ such that $u_k$ tends to $-\infty$
(resp. $+\infty$) on $[Q_i^k,P_{i+1}]$ (resp. on $[P_i,Q_i^k]$) and
$u_k(P_0)=0$. Following  Remark \ref{aa} in section \ref{seqmse}, we
shall prove 
  that this sequence $(u_k)$ of 
solutions of \eqref{MSE} has no line of divergence, then the limit $u$
of $(u_k)$ will be our solution. We shall make discussions that are
similar to  the ones made by H.~Jenkins and J.~Serrin. We note
$\dd\Psi_{u_k}=\dd\Psi_k$. We recall that, if $T$ is a segment included
in a line of divergence, $\left|\dis\int_T\dd\Psi_k\right|$ converge
to the length of $T$ for a subsequence.

Suppose there exists a line of divergence $L$. We first prove that $L$ 
can not have an end point in the interior of a $L_i^+$ or a
$L_i^-$. Suppose that $L$ has an end point $D$ in $L_i^-$ (the same
argument works for $L_i^+$). Let $A$ be a point in $L\cap S_{i-1}$, we
orient $L$ by $\overrightarrow{AD}$, we suppose that the limit normal
along $L$ points on the right-hand side of $L$. We chose a point $B$ of
$L_i^-$ on the right-hand side of $D$. Because of the triangle inequality,
there exists a point $C$ on $[A,D]$ such that
$|AC|+|DB|>|CD|+|BA|$. for $k$ big enough we have $A$ and $C$ in
$\Ome_k$ we then put $D_k=[A,D]\cap[P_i,Q_{i-1}^k]$ and
$B_k=[A,B]\cap[P_i,Q_{i-1}^k]$ (see Figure \ref{plusinf}). 

\begin{figure}[h]
\begin{center}
\resizebox{1\linewidth}{!}{\input{figure1.pstex_t}}
%\scalebox{0.8}{\input{figure1.pstex_t}}
\caption{\label{plusinf}}
\end{center}
\end{figure}

Let $T_k$ be the triangle $AD_kB_k$ with this orientation. We then
have:

\begin{equation}
\begin{split}
0=\int_{T_k}\dd\Psi_k & \ge \int_{[A,C]}\dd\Psi_k
-|CD_k|+|D_kB_k|-|B_kA|\\
&\ge \int_{[A,C]}\dd\Psi_k -|CD|+|D_kB_k|-|BA| \\
\end{split}
\end{equation}

But $|D_kB_k|\longrightarrow|DB|$ and $\dis \int_{[A,C]}\dd\Psi_k
\longrightarrow |AC|$ for the subsequence that makes $L$ appear; this 
gives us a contradiction.

We have now only a finite number of possibilities for a line of
divergence. If it has an end point, it must be one $P_i$. By
construction, we have $\dis\int_{[P_i,P_{i+1}]}\dd\Psi_k=0$ (because
the integral of $\dd\Psi_k$ along the triangle $P_iQ_i^kP_{i+1}$ is
zero and we know $\dd\Psi_k$ along $[P_iQ_i^k]$ and $[Q_i^kP_{i+1}]$
by Lemma \ref{calc}) so if
$\Gamma$ is a curve joining $P_i$ to $P_j$ we have
$\dis\int_\Gamma\dd\Psi_k=0$. Then, by passing to the limit, if
$\Gamma$ is a line of divergence, we obtain
$|P_iP_j|=0$ which is not possible. This proves that a line of
divergence 
has at most one end point. Suppose that a line of divergence $L$ has
no end point, we are in the situation of Figure \ref{noextrem}. Let
$A$ and $B$ 
be point on $L$ as in Figure \ref{noextrem} such that
$|AB|>|P_iP_{i+1}|$. We note $D$  
(resp. $C$) the projection of $A$ (resp. $B$) on $L_{i+1}^-$. For $k$
big enough we note $C_k=[B,C]\cap[P_{i+1},Q_i^k]$ and $D_k= [A,D]\cap
[P_{i+1},Q_i^k]$. We then have:

\begin{equation}
\begin{split}
0=\int_{ABC_kD_k}\dd\Psi_k &\ge \int_{[A,B]}\dd\Psi_k -|BC_k|
+\int_{[C_k,D_k]}\dd\Psi_k -|D_kA|\\
&\ge \int_{[A,B]}\dd\Psi_k -2|P_i,P_{i+1}|+|C_kD_k|
\end{split}
\end{equation} 

We have $|C_kD_k|\longrightarrow|CD|=|AB|$ and $\dis\int_{[A,B]}\dd\Psi_k
\longrightarrow |AB|$ for a subsequence, so we get a contradiction and
a line of divergence must have one end point.

\begin{figure}[h]
\begin{center}
\resizebox{1\linewidth}{!}{\input{figure4.pstex_t}}
\caption{\label{noextrem}}
\end{center}
\end{figure}

Let $L$ be a line of divergence, we know that we are in the case where
$L$ has $P_i$ as end point and goes
to infinity in one $S_j$. By what we have done just above, we have only
one possiblity for the limit normal: we are in the same situation as in the
semi-strip $S_j$ in Figure \ref{noextrem}. Then, by changing $L$ if
necessary, 
we can suppose that the part of $L$ in $S_j$ is parametrized by
$\{A\}\times\R_+$ with $A\in]P_j,P_{j+1}[$ and the domain
$\widetilde{\Ome}$ parametrized by $[A,P_{j+1}]\times\R_+$ is in
$\boB(u_k)$. Let $\theta$ be an extraction that makes $L$ appear,
since $\widetilde{\Ome}\subset\boB(u_k)$, there exists an extraction
$\theta'$ such that $u_{\theta'(k)}$ is a subsequence of
$u_{\theta(k)}$ and $u_{\theta'(k)}-u_{\theta'(k)}(K)$ (where
$K\in\widetilde{\Ome}$) converges to $v$ a solution of \eqref{MSE} on
$\widetilde{\Ome}$.

We shall now prove that $v$ tends to $-\infty$ on $L_{j+1}^-$ and
$+\infty$ on $L$. Let $B\in[A,P_{j+1}]$ and C and D be the points which are
respectively parametrized by $(B,c)$ and $(B,d)$ ($c<d$). We note:

\begin{gather*}
\textrm{$E$ the projection of $D$ on $L_{j+1}^-$}\\
\textrm{$F$ the projection of $C$ on $L_{j+1}^-$}\\
\textrm{$G$ the projection of $D$ on $L$}\\
\textrm{$H$ the projection of $C$ on $L$}
\end{gather*}
We note also, when $k$ is big enough, $E_k= [D,E]\cap [P_{j+1},Q_j^k]$
and $F_k= [C,F]\cap [P_{j+1},Q_j^k]$. Because $\dd\Psi_k$ is closed we
have:
\begin{gather}
\left|\int_{[C,D]}\dd\Psi_{\theta'(k)}- |F_{\theta'(k)} E_{\theta'(k)}|
  \right| \le 2|BP_{j+1}|\\
\left|\int_{[C,D]}\dd\Psi_{\theta'(k)}- \int_{[H,G]}\dd\Psi_{\theta'(k)}
  \right| \le 2|BA|
\end{gather} 
Thus, in letting $k$ tends to infinity, we obtain:

\begin{gather}
\left|\int_{[C,D]}\dd\Psi_v-|CD|\right| \le 2|BP_{j+1}|\\
\left|\int_{[C,D]}\dd\Psi_v-|CD|\right| \le 2|BA|
\end{gather}

So we can calculate $\dd\Psi_v$ on $L$ and $L_{j+1}^-$, we remark that
$\dd\Psi_v$ has the same behaviour as if $v$ assumes the boundary values
$+\infty$ on $L$ and $-\infty$ on $L_{j+1}^-$. We prove that this is, in fact,
the case. We consider
now two points $A_1$ and $A_2$ on $L\cap\widetilde{\Omega}$ and two
points $A_3$ and $A_4$ on $L_{j+1}^-$. There exists a solution $v'$ of
\eqref{MSE} on the domain bounded by the polygon $A_1A_2A_3A_4$ such
that $v'=v$ on $[A_1,A_4]$ and $[A_2,A_3]$, $v'$ tends to $+\infty$ on
$[A_1,A_2]$ and tends to $-\infty$ on $[A_3,A_4]$. Since we know the
value of $\dd\Psi_v$ on $[A_1,A_2]$ and $[A_3,A_4]$, the uniqueness part
of the proof of Theorem \ref{thmjs} proves that $v=v'$. We then have
proved that $v$ tends to $+\infty$ on $L$ and $-\infty$ on $L_{j+1}^-$. 

We shall now get a contradiction to the existence of the line of divergence
$L$. We have $\dis\int_{[A,P_{j+1}]}\dd\Psi_v= 
\int_{\widetilde{\Ome}\cap I_{j,l}} \dd\Psi_v$, then, by \eqref{bb}
and letting $l$ tends to infinity, we get
$\dis\int_{[A,P_{j+1}]}\dd\Psi_v=0$. If we follow $L$ between $P_i$
and $A$ and the segment $[A,P_{j+1}]$, we get a path joining $P_i$ to
$P_{j+1}$. Then we have:

\begin{equation}
0=\int_{[P_i,A]}\dd\Psi_{\theta'(k)}+\int_{[A,P_{j+1}]}\dd\Psi_{\theta'(k)}
\end{equation}

Let $k$ tend to infinity, we get $0=|P_iA|+\dis \int_{[A,P_{j+1}]}
\dd\Psi_v= |P_iA|$; but $P_i\notin[P_j,P_{j+1}]$, this is our
contradiction. 

We then have prove that $\boB(u_k)=\Ome(\boP)$, as $u_k(P_0)=0$ for
all $k$ there exists a subsequence $u_{k'}$ which converges to a
solution $u$ of \eqref{MSE}. The same arguments that we used just
above for $v$ prove
that $u$ tends to $+\infty$ (resp. $-\infty$) on $L_i^-$
(resp. $L_i^+$); we have then established Theorem \ref{moi}. \qed

\bigskip

We are then able to build the solution to the Plateau problem at
infinity. Let $V$ be a polygon and $(\boP,\phi)$ a polygonal disk bounded by
$V$. We consider the solution $u$ of the Dirichlet problen given by
Theorem \ref{moi}. We note $P_1,\dots,P_r$ the vertices of $V$, we
consider $\Psi_u$ normalized by $\Psi_u(P_1)=0$, from the proof above, we have
$\Psi_u(P_i)=0$ for all $i$. Then on $L_i^+$ and $L_i^-$ we have
$\Psi_u(Q)=|QP_i|$; since $\Psi_u$ is $1$-Lipschitz continuous, we
have $\Psi_u(Q)\ge 0$ for all $Q\in S_i$. Suppose that
$\{Q\in\Ome(\boP)|\ \Psi_u(Q)\le 0\}$ is not reduced to
$\{P_1,\dots,P_r\}$ then there exists a point in the interior of
$\boP$ such that $\Psi_u$ is minimal at this point. But $\Psi_u$
corresponds to $x_3^*$, the third coordinate on $M^*$ the conjugate
surface to the graph of $u$; since $x_3^*$ is harmonic on $M^*$, it
can not have a minimum in the interior of $M^*$. We then have proved
that $\Psi_u>0$ in the interior of $\Ome(\boP)$.

By Theorem \ref{regbound}, the boundary of the graph of $u$ is composed of the
$r$ vertical lines over the points $\phi(P_i)$. Let $M$ be the graph
of $u$ with 
these $r$ vertical lines. We consider $M^*$ the conjugate 
surface to $M$. The boundary of $M^*$ is composed of $r$ horizontal
planar geodesic curves, since $\Psi_u(P_i)=0$ for all $i$ the $r$
curves are 
all in the plane $\{x_3=0\}$. Finally, we consider $\Sigma$ the union
of $M^*$ and of its symmetry by $\{x_3=0\}$. The surface $\Sigma$ is a
regular minimal surface, it is complete and its flux polygon is $V$ by
construction. By construction, we know also that $\Sigma$ is strongly
symmetric with respect to $\{x_3=0\}$.

The last thing we have to prove
about $\Sigma$ for being sure that it is the solution of the Plateau
problem at infinity is that it has finite total curvature.

We know (see \cite{Os}) that there exists a constant $c$ such that if
$u$ is a solution of \eqref{MSE} on a domain $D$ and $A\in D$, $M$ is
the graph of $u$ and $d$ is the distance along $M$ of the point in $M$
over $A$ to the boundary of $S$ then the curvature $K$ of $M$ at the point
over $A$ is bounded by $\dis\frac{c}{d^2 W^2(A)}$.

Let us consider a half-strip $S=[0,a]\times \R_+$ and $u$ a solution of
\eqref{MSE} on $S$ such that $u$ takes the value $+\infty$
(resp. $-\infty$) on $\{0\}\times \R_+^*$ (resp. $\{a\}\times \R_+^*$).
The boundary of the graph of $u$ is over $[0,a]\times\{0\}$. We then
have $\dis K(x,y)\le \frac{c}{x^2W^2(x,y)}$. We consider the part
$S'\subset{S}$ such that $x\ge x_0>0$. For a domain $D$ we note $K(D)$ the
total curvature of the graph over $D$. We then have:

\begin{equation*}
\begin{split}
K(S')&=\int_{S'}K(x,y)W(x,y)\dd x\dd y\\
&\le\int_{S'}\frac{c}{x^2W(x,y)}\dd x\dd y\\
&\le\int_{x_0}^{+\infty}\frac{ca}{x^2}\dd x=\frac{ca}{x_0}<+\infty
\end{split}
\end{equation*}

We now use arguments that are similar to the first part of the proof of
Theorem \ref{regbound}. We consider, for $\alpha\in]0,\frac{\pi}{2}[$,

$S(\alpha)=\{(x,y)\in S|\ a-y\le x\tan\alpha \}$. Let us take $\alpha$ such
that $\tan 
\alpha<\frac{1}{8}$, then Lemma 1 of \cite{JS} proves that for every $(x,y)\in
S(\alpha)$ $q(x,y)<0$. We note $L(\alpha)$ the segment in $S(\alpha)$ such
that $a-y=x\tan\alpha$. By Lemma \ref{psherck}, $u$ is lower-bounded by $m_1$
on $L(\alpha)$ and upper bounded by $m_2$ on the part of $L$ such that
$x<x_0<a$. We then define $\Theta:(x,y)\mapsto\left(x,u(x,y)\right)$. $\Theta$
is a diffeomorphism of $S(\alpha)$ into its image $\im\Theta$. We define
$\chi=\Theta^{-1}$ then $\chi_{(2)}$ is a solution of \ref{MSE} on $\im
\Theta$. We observe that $\chi$ extends smoothly to $\Theta(L)$. We have
$\chi_{(2)}(x,z)$ tends to $a$ as $x$ tends to $0$ so we can extend
$\chi_{(2)}$ by symmetry to $\im\Theta\cup\{(x,z)\in\R^2|
(-x,z)\in\im\Theta\}$. To compute the total curvature of the graph of $u$
over $D(\alpha)\cap\{x<x_0\}$, we use its parametrization as a graph over
$\im\Theta\cap\{0\le x\le x_0\}$:

\begin{equation*}
\begin{split}
K(D(\alpha)\cap\{x<x_0\})&= K(\im\Theta\cap\{0\le x\le x_0\})\\
&=\int_{\im\Theta\cap\{0\le x\le x_0\}\cap\{z\le m_1-1\}}K(x,z)W(x,z)\dd x\dd
z\\ 
&\quad+\int_{\im\Theta\cap\{0\le x\le x_0\}\cap\{z\ge m_1-1\}}K(x,z)W(x,z)\dd
x\dd z\\ 
&\le \int_{\im\Theta\cap\{0\le x\le x_0\}\cap\{z\le
  m_1-1\}}\frac{c}{(z-m_1)^2W(x,z)}\dd x\dd z\\
&\quad+\int_{\im\Theta\cap\{0\le x\le x_0\}\cap\{z\ge m_1-1\}}K(x,z)W(x,z)\dd
x\dd z\\    
&\le \int_1^{+\infty}\frac{cx_0}{z^2}\dd z+C\\
\intertext{because $\overline{\im\Theta}\cap\{0\le x\le x_0\}\cap\{z\ge
  m_1-1\}$ is compact} 
&<+\infty
\end{split}
\end{equation*}

We can do the same work for $\{(x,y)\in S| \ y\le x\tan \alpha\}$.

We then control the curvature on each semi-strip $S_i$. There is a last
part in $\Ome(\boP)$. This part is compact and by Lemma \ref{psherck}, $u$ is
bounded on this part; besides the graph is regular at the boundary. So the
graph above this last part is a compact part of the whole graph then it has
finite total curvature. We then have proved that the graph $M$ has
finite total curvature. Since $M^*$ is isometric to $M$, it has finite
total curvature and then $\Sigma$ has finite total curvature because
it is twice as many as the one of $M^*$.

%%% Local Variables: 
%%% mode: latex
%%% TeX-master: t
%%% End: 

%% file: figure1.pstex_t
\begin{picture}(0,0)%
\epsfig{file=figure1.pstex}%
\end{picture}%
\setlength{\unitlength}{3947sp}%
\begingroup\makeatletter\ifx\SetFigFont\undefined%
\gdef\SetFigFont#1#2#3#4#5{%
  \reset@font\fontsize{#1}{#2pt}%
  \fontfamily{#3}\fontseries{#4}\fontshape{#5}%
  \selectfont}%
\fi\endgroup%
\begin{picture}(10662,4122)(826,-6373)
\put(2926,-2461){\makebox(0,0)[lb]{\smash{\SetFigFont{14}{16.8}{\rmdefault}{\mddefault}{\updefault}$L_i^-$}}}
\put(2776,-6286){\makebox(0,0)[lb]{\smash{\SetFigFont{14}{16.8}{\rmdefault}{\mddefault}{\updefault}$L$}}}
\put(3451,-4861){\makebox(0,0)[lb]{\smash{\SetFigFont{14}{16.8}{\rmdefault}{\mddefault}{\updefault}$A$}}}
\put(5851,-2536){\makebox(0,0)[lb]{\smash{\SetFigFont{14}{16.8}{\rmdefault}{\mddefault}{\updefault}$D$}}}
\put(9676,-2461){\makebox(0,0)[lb]{\smash{\SetFigFont{14}{16.8}{\rmdefault}{\mddefault}{\updefault}$B$}}}
\put(5251,-2986){\makebox(0,0)[lb]{\smash{\SetFigFont{14}{16.8}{\rmdefault}{\mddefault}{\updefault}$D_k$}}}
\put(7801,-3661){\makebox(0,0)[lb]{\smash{\SetFigFont{14}{16.8}{\rmdefault}{\mddefault}{\updefault}$B_k$}}}
\put(5101,-3661){\makebox(0,0)[lb]{\smash{\SetFigFont{14}{16.8}{\rmdefault}{\mddefault}{\updefault}$C$}}}
\put(826,-2461){\makebox(0,0)[lb]{\smash{\SetFigFont{14}{16.8}{\rmdefault}{\mddefault}{\updefault}$P_i$}}}
\end{picture}

%% file: figure4.pstex_t
\begin{picture}(0,0)%
\epsfig{file=figure4.pstex}%
\end{picture}%
\setlength{\unitlength}{4144sp}%
\begingroup\makeatletter\ifx\SetFigFont\undefined%
\gdef\SetFigFont#1#2#3#4#5{%
  \reset@font\fontsize{#1}{#2pt}%
  \fontfamily{#3}\fontseries{#4}\fontshape{#5}%
  \selectfont}%
\fi\endgroup%
\begin{picture}(11949,4749)(664,-5473)
\put(6931,-3931){\makebox(0,0)[lb]{\smash{\SetFigFont{12}{14.4}{\rmdefault}{\mddefault}{\updefault}$P_i$}}}
\put(6256,-2176){\makebox(0,0)[lb]{\smash{\SetFigFont{12}{14.4}{\rmdefault}{\mddefault}{\updefault}$P_{i+1}$}}}
\put(2206,-4471){\makebox(0,0)[lb]{\smash{\SetFigFont{12}{14.4}{\rmdefault}{\mddefault}{\updefault}$P_{j+1}$}}}
\put(2341,-2311){\makebox(0,0)[lb]{\smash{\SetFigFont{12}{14.4}{\rmdefault}{\mddefault}{\updefault}$P_j$}}}
\put(5491,-3166){\makebox(0,0)[lb]{\smash{\SetFigFont{12}{14.4}{\rmdefault}{\mddefault}{\updefault}L}}}
\put(11521,-3121){\makebox(0,0)[lb]{\smash{\SetFigFont{12}{14.4}{\rmdefault}{\mddefault}{\updefault}A}}}
\put(7741,-2086){\makebox(0,0)[lb]{\smash{\SetFigFont{12}{14.4}{\rmdefault}{\mddefault}{\updefault}C}}}
\put(7786,-3076){\makebox(0,0)[lb]{\smash{\SetFigFont{12}{14.4}{\rmdefault}{\mddefault}{\updefault}B}}}
\put(11566,-2086){\makebox(0,0)[lb]{\smash{\SetFigFont{12}{14.4}{\rmdefault}{\mddefault}{\updefault}D}}}
\end{picture}